
\input amstex.tex
\documentstyle{amsppt}
\magnification1200
\hsize=12.5cm
\vsize=18cm
\hoffset=1cm
\voffset=2cm

\def\DJ{\leavevmode\setbox0=\hbox{D}\kern0pt\rlap
{\kern.04em\raise.188\ht0\hbox{-}}D}
\def\dj{\leavevmode
 \setbox0=\hbox{d}\kern0pt\rlap{\kern.215em\raise.46\ht0\hbox{-}}d}

\def\txt#1{{\textstyle{#1}}}
\baselineskip=13pt
\def\hf{{\textstyle{1\over2}}}
\def\a{\alpha}
\def\d{{\,\roman d}}
\def\e{\varepsilon}

\def\s{\sigma}
\def\t{\theta}
\def\={\;=\;}

\def\zt{\zeta(\hf+it)}

\def\R{\Re{\roman e}\,} 
\def\z{\zeta}

 \def\t{\theta}
\def\hf{{\textstyle{1\over2}}}
\def\txt#1{{\textstyle{#1}}}

\font\tenmsb=msbm10
\font\sevenmsb=msbm7
\font\fivemsb=msbm5
\newfam\msbfam
\textfont\msbfam=\tenmsb
\scriptfont\msbfam=\sevenmsb
\scriptscriptfont\msbfam=\fivemsb
\def\Bbb#1{{\fam\msbfam #1}}

\def \NN {\Bbb N}

\def \RR {\Bbb R}

\font\ff=cmr8
\def\txt#1{{\textstyle{#1}}}
\baselineskip=13pt

\font\teneufm=eufm10
\font\seveneufm=eufm7
\font\fiveeufm=eufm5
\newfam\eufmfam
\textfont\eufmfam=\teneufm
\scriptfont\eufmfam=\seveneufm
\scriptscriptfont\eufmfam=\fiveeufm
\def\mathfrak#1{{\fam\eufmfam\relax#1}}

\font\tenmsb=msbm10
\font\sevenmsb=msbm7
\font\fivemsb=msbm5
\newfam\msbfam
     \textfont\msbfam=\tenmsb
      \scriptfont\msbfam=\sevenmsb
      \scriptscriptfont\msbfam=\fivemsb
\def\Bbb#1{{\fam\msbfam #1}}

\def \NN {\Bbb N}

\def \RR {\Bbb R}

  \def\rightheadline{{\hfil{\ff
Moments of the zeta-function in short intervals}\hfil\tenrm\folio}}

  \def\leftheadline{{\tenrm\folio\hfil{\ff
 Aleksandar Ivi\'c }\hfil}}
  \def\emptyheadline{\hfil}
  \headline{\ifnum\pageno=1 \emptyheadline\else
  \ifodd\pageno \rightheadline \else \leftheadline\fi\fi}

\font\ff=cmr8
\font\teneufm=eufm10
\font\seveneufm=eufm7
\font\fiveeufm=eufm5
\newfam\eufmfam
\textfont\eufmfam=\teneufm
\scriptfont\eufmfam=\seveneufm
\scriptscriptfont\eufmfam=\fiveeufm
\def\mathfrak#1{{\fam\eufmfam\relax#1}}

\font\tenmsb=msbm10
\font\sevenmsb=msbm7
\font\fivemsb=msbm5
\newfam\msbfam
\textfont\msbfam=\tenmsb
\scriptfont\msbfam=\sevenmsb
\scriptscriptfont\msbfam=\fivemsb
\def\Bbb#1{{\fam\msbfam #1}}

\def \NN {\Bbb N}

\def \RR {\Bbb R}

\def\a{\alpha}
 \def\e{\varepsilon}
 \def\d{\,{\roman d}}
\topmatter
\title
On the moments of the Riemann zeta-function in short  intervals
\endtitle
\author
 Aleksandar Ivi\'c
\endauthor
\address
Katedra Matematike RGF-a,
Universitet u Beogradu,  \DJ u\v sina 7,
11000 Beograd, Serbia.\bigskip
\endaddress
\keywords  power moments, Riemann zeta-function, Riemann Hypothesis,
upper bounds
\endkeywords
\subjclass 11 M 06
\endsubjclass
\email
{\tt
ivic\@rgf.bg.ac.rs,\enskip aivic\@matf.bg.ac.rs}
\endemail
\dedicatory
See Hardy-Ramanujan J. {\bf32}(2009), 4-23, www.nias.res.in/hrj/contentsvol32.htm. 
\enddedicatory
\abstract Assuming the Riemann Hypothesis it is proved that,
for fixed $k>0$ and $H = T^\t$ with fixed $0<\t \le 1$,
$$
\int_T^{T+H}|\zt|^{2k}\d t \ll H(\log T)^{k^2(1+O(1/\log_3T))},
$$
where $\log_jT = \log(\log_{j-1}T)$. The proof is based on the method of
K. Soundararajan [8] for counting the occurrence of large values
of $\log|\zt|$, who proved that
$$
\int_0^{T}|\zt|^{2k}\d t \ll_\e T(\log T)^{k^2+\e}.
$$

\endabstract
\endtopmatter
\document
\head
1. Introduction
\endhead
Power moments of $|\zt|$ are a central problem in the theory of the
Riemann zeta-function $\z(s) = \sum\limits_{n=1}^\infty n^{-s}\;(\s = \R s >1)$
and a vast literature
exists on this subject (see e.g., the monographs [3], [4], [6] and [9]). If
$$
I_k(T,H) := \int_T^{T+H}|\zt|^{2k}\d t\qquad(k >0,\; 1 \ll H \le T),\leqno(1.1)
$$
then naturally one seeks asymptotic formulas for $I_k(T,H)$.
It is only in the cases $k = 1$ and $k=2$ that we have precise expressions for
$I_k(T,T)$, due to the well-known works of F.V. Atkinson [1] and Y. Motohashi
(see e.g. [6]), respectively. Although with the use of methods relating to random
matrix theory (see J.B. Conrey et al. [2]) it is possible to make plausible
conjectures for the asymptotic formulas for $I_k(T,T)$ when $k\in\NN$ is fixed,
no one has proved yet such an asymptotic formula for $k\ge3$, even assuming the
Riemann Hypothesis (RH, that all complex zeros of $\z(s)$ have real parts
equal to 1/2).
Unconditional lower bounds for $I_k(T,H)$ and similar expressions involving the
derivatives of $\z(s)$ have been obtained in a series of papers by R. Balasubramanian
and K. Ramachandra. These results, which are of a general nature and involve various
convexity techniques, are expounded in Ramachandra's monograph [7].
For example, one has unconditionally
$$
I_k(T,H) \;\gg \;H(\log H)^{k^2}\qquad(k\in\NN,\; \log_2T \ll H \le T)\leqno(1.2)
$$
when $k$ is fixed and $\log_jT = \log(\log_{j-1}T)$ is the $j$-th iteration of
the natural logarithm. Under the RH it is known that (1.2) holds for any fixed $k>0$.

\medskip
Furthermore a classical result of  J.E. Littlewood states (see e.g., [9] for a proof)
that, under the RH,
$$
|\zt| \ll \exp\left(C\,{\log t\over\log_2t}\right)\qquad\Bigl(C>0\Bigr),\leqno(1.3)
$$
which can be used to provide a trivial upper bound for $I_k(T,H)$.
However, recently K. Soundararajan [8] complemented (1.2) in the case $H=T$ by
obtaining, under the RH, the non-trivial upper bound
$$
\int_0^T|\zt|^{2k}\d t \; \ll_{\e}\; T(\log T)^{k^2+\e},\leqno(1.4)
$$
which is valid for any fixed $k>0$ and any given $\e>0$.
This result, apart from `$\e$', is therefore best possible.
His method of proof is based on a large values
estimate for $\log|\zt|$, which gives
as a corollary (under the RH) the bound (1.3) with the explicit constant $C = 3/8$.

\medskip
The aim of this note is twofold. The main objective is to generalize
(under the RH) (1.4) to
upper bounds for $I_k(T,H)$. The second aim is to replace `$\e$' by an
explicit function of $T$ which is $o(1)$ as $T\to\infty$. The result is
contained in

\medskip
THEOREM 1. {\it Let $H = T^\t$ where $0 < \t \le 1$ is a fixed number, and
let $k$ be a fixed positive number.Then, under the RH,   we have}
$$
I_k(T,H) \;=\; \int_T^{T+H}|\zt|^{2k}\d t
\;\ll\; H(\log T)^{k^2\bigl(1+O(1/\log_3T)\bigr)}.
\leqno(1.5)
$$

\medskip
Note that, since $H = T^\t$ and $\t$ is fixed, the right-hand side of (1.5)
is
$$
\ll H(\log H)^{k^2\bigl(1+O(1/\log_3T)\bigr)},
$$
which is more in tune with the lower bound in (1.2), although it does not seem
possible to reach the range $\log_2T \ll H\le T$ by the present method
or to remove $O(1/\log_3T)$ from the exponent in (1.5). As already noted,
upper bounds of the form
$I_k(T,H) \ll H(\log T)^{k^2}$ can be derived unconditionally in the cases
$k= 1$ and $k =2$. They are known to hold for $\t>1/3$ (and even for some slightly
smaller values of $\t$) when $k=1$, and  for $\t>2/3$ when $k=2$. In the case
when $k=1/2$ it is known (see K. Ramachandra [7]) that this bound holds unconditionally
when $\t>1/2$  and for $\t>1/4$ under the RH. No other results of this type seem to
be known for other values of $\t$.

\smallskip

Theorem 1 will be deduced from a large values estimate for $\log|\zt|$,
based on Soundararajan's method [8]. This is
\medskip
THEOREM 2. {\it Let $H = T^\t$ where $0 < \t \le 1$ is a fixed number, and
let $\mu(T,H,V)$ denote the measure of points $t$ from $[T,\,T+H]$ such that
$$
\log|\zt| \ge V,\quad 10\sqrt{\log_2T} \le V \le {3\log 2T\over 8\log_2(2T)}.
\leqno(1.6)
$$
Then, under the RH, for $\;10\sqrt{\log_2T} \le V \le\log_2T$  we have
$$
\mu(T,H,V) \;\ll\; H{V\over\sqrt{\log_2T}}\exp\left(-{V^2\over\log_2T}\Bigl(1-
{7\over2\t\log_3T}\Bigr)\right),
\leqno(1.7)
$$
for $\;\log_2T \le V \le\hf\t\log_2T\log_3T$ we have
$$
\mu(T,H,V) \;\ll\; H\exp\left(- {V^2\over\log_2T}
\Bigl(1 - {7V\over4\t\log_2T\log_3T}\Bigr)^2\,\right),\leqno(1.8)
$$
and for $\;\hf\t\log_2T\log_3T\le V\le {3\log 2T\over 8\log_2(2T)}$ we have}
$$
\mu(T,H,V) \;\ll\; H\exp(-\txt{1\over20}\t V\log V).\leqno(1.9)
$$

\bigskip
\head
2. The necessary Lemmas
\endhead
In this section we shall state the necessary lemmas for the
proof of Theorem 1 and Theorem 2.

\medskip
LEMMA 1. {\it Assume the RH. Let $T\le t\le 2T, T\ge T_0, 2 \le x \le T^2$.
If $\lambda_0 = 0.4912\ldots\,$ denotes the unique positive real
number satisfying $\,{\roman e}^{-\lambda_0} = \lambda_0 + \hf\lambda_0^2$, then
for $\lambda \ge \lambda_0$ we have}
$$
\log|\zt| \le \R\sum_{2\le n\le x}{\Lambda(n)\over n^{{1\over2}+{\lambda\over \log x}
+ it}\log n}{\log(x/n)\over \log x} + {(1+\lambda)\over2}{\log T\over \log x}
+ O\Bigl({1\over\log x}\Bigr).\leqno(2.1)
$$

\medskip Lemma 1 is due to K. Soundararajan [8]. It is based on Selberg's
classical method (see e.g., E.C. Tichmarsh [9, Th. 14.20]) of the use of
an explicit expression for $\z'(s)/\z(s)$ by means of a sum containing the
familiar von Mangoldt function $\Lambda(n)$ (equal to $\log p$ if $n = p^\a$,
where $p$ denotes primes, and $\Lambda(n) = 0$ otherwise). Soundararajan's
main innovation is the observation that
$$
F(s) := \R \sum_\rho{1\over s-\rho} = \sum_\rho {\s-\hf\over (\s-\hf)^2+ (t-\gamma)^2}
\ge 0\qquad(s = \s+it)
$$
if $\s \ge \hf$, since all the complex zeros $\rho$ of $\z(s)$ are (this is  the RH) of
the form $\rho = \hf + i\gamma\;(\gamma \in\RR)$. Note that $F(s)$ appears
in the classical expression
(this is unconditional, valid when $s=\s+it$, $t$ is not an ordinate of any $\rho$ and
$T\le t\le 2T$)
$$
-\R {\z'(s)\over\z(s) }= \hf\log T - F(s) + O(1).
$$
This is integrated over $\s$ from $\hf$ to $\s_0\,(>\hf)$ to produce ($s_0 = \s_0+it$)
$$
\log|\zt| - \log|\z(s_0)| \le (\s_0 - \hf)\left(\hf\log T - \hf F(s_0) + O(1)\right),
$$
at which point an expression similar to [9, Th. 14.20] for $\log|\z(s_0)|$ is used,
and the non-negativity of $F(s)$ can be put to advantage.

\medskip
LEMMA 2. {\it Assume the RH. If $T \le t \le 2T, 2 \le x \le T^2, \s \ge \hf$, then}
$$
\sum_{2\le n \le x, n\ne p}{\Lambda(n)\over n^{\s
+ it}\log n}{\log(x/n)\over \log x} \;\ll\; \log_3T.\leqno(2.2)
$$

\medskip This is Lemma 2 from Soundararajan [8], where a brief sketch of the proof
is indicated. The details are to be found in the work of M.B. Milinovich [5].

\medskip
LEMMA 3. {\it Let $\,2 \le x \le T,\, T \ge T_0$. Let $1 \ll H \le T$ and $r\in \NN$
satisfy $x^r \le H$. For any complex numbers $a(p)$ we have}
$$
\int_T^{T+H}\Bigl|\sum_{p\le x}{a(p)\over p^{{1\over2}+it}}\Bigr|^{2r}\d t
\;\ll\;Hr!{\Bigl(\sum_{p\le x}{|a(p)|^2\over p}\Bigr)}^r.\leqno(2.3)
$$

\medskip
{\bf Proof}. This lemma is unconditional and
is a generalization of Lemma 3 of Soundararajan [8] (when $H=T$),
but he had the more stringent condition $x^r \le T/\log T$ (we have changed
his notation from $k$ to $r$ to avoid confusion with $k$ in Theorem 1). Write
$$
\Bigl(\sum_{p\le x}{a(p)\over p^{{1\over2}+it}}\Bigr)^r
= \sum_{n\le x^r}{a_{r,x}(n)\over n^{{1\over2}+it}},
$$
where $a_{r,x}(n)=0$ unless $n$ is a product of $r$ primes factors,
each of which is $\le x$.
By the mean value theorem for Dirichlet polynomials (see e.g., [3, Chapter 4])
the left-hand side of (2.3) is equal to
$$
H\sum_{n\le x^r}{|a_{r,x}(n)|^2\over n}
+ O\Bigl(\sum_{n\le x^r}|a_{r,x}(n)|^2\Bigr) \ll (H + x^r)
\sum_{n\le x^r}{|a_{r,x}(n)|^2\over n}.\leqno(2.4)
$$
But, as shown in detail in [8], it is not difficult  to see that
$$
\sum_{n\le x^r}{|a_{r,x}(n)|^2\over n} \;\ll\;
r!{\Bigl(\sum_{p\le x}{|a(p)|^2\over p}\Bigr)}^r,
$$
hence (2.3) follows from (2.4), since our assumption is that $x^r \le H$.

\head
3. Proof of  Theorem 1
\endhead
The contribution of $t$ satisfying $\log|\zt| \le \hf k\log_2T$
to $I_k(T,H)$ is trivially
$$
\le \;H\Bigl\{(\log T)^{k/2}\Bigr\}^{2k} = H(\log T)^{k^2}.\leqno(3.1)
$$
Likewise the bound (1.5) holds, by (1.7) and (1.8), for the
contribution of $t$ satisfying $\log|\zt| \ge 10k\log_2T$. Thus we can
consider only the range
$$
V + {j-1\over\log_3T} \;\le\;\log|\zt|\;\le\; V + {j\over\log_3T},\leqno(3.2)
$$
where $1 \le j \ll \log_3T, V = 2^{\ell-{1\over2}}k\log_3T,
\, 1 \le \ell \le {3\over2} + \bigl[{\log10\over\log2}\bigr]$. If we set
$$
U = U(V,j;T)\;:=\; V + {j-1\over\log_3T},\leqno(3.3)
$$
then we have
$$
I_k(T,H) \ll H(\log T)^{k^2} + \log_3T\max_U\mu(T,H,U)
\exp\Bigl(2k(U + 1/\log_3T)\Bigr),
\leqno(3.4)
$$
where $\mu(T,H,U)$ is the measure of $t\in [T,T+H]$ for which $\log|\zt| \ge U$,
and the maximum is over $U$ satisfying (3.2)--(3.3).
If we use (1.7) and (1.8) of Theorem 2, then in the relevant range for $U$ we obtain
$$\eqalign{
\mu(T,H,U)\exp\bigl(2k(U + 1/\log_3T)\bigr) &
\ll H\log_2T \exp\Bigl(2kU - U^2G(T)\Bigr),\cr
G(T) &\;:=\; {1\over\log_2T}\Bigl(1 + O\bigl({1\over\log_3T}\bigr)\Bigr).\cr}
$$
Since $2kU - U^2G(T)$ attains its maximal value at $U = k/G(T)$, we have
$$\eqalign{&
\mu(T,H,U)\exp\bigl(2k(U + 1/\log_3T)\bigr) \cr&\ll
H\log_2T\,\exp\Bigl(k^2\log_2T(1+O(1/\log_3T)\Bigr)\cr&
= H(\log T)^{k^2(1+O(1/\log_3T))},\cr}
$$
so that (3.4) yields then (1.5) of Theorem 1.

\break

\head
4. Proof of  Theorem 2
\endhead
We let
$$
x\= H^{A/V},\quad z\= x^{1/\log_2T},\quad A = A(T,V) \;(\;\ge 1),\leqno(4.1)
$$
where $A$ will be suitably chosen below. We follow the method of proof of
[8] and accordingly consider three cases.

\medskip
{\bf Case 1}. When $10\sqrt{\log_2T} \le V \le \log_2T$, we take $A = \hf\log_3T$.

\medskip
{\bf Case 2}. When $\log_2T \le V \le \hf \t\log_2T\log_3T$, we take
$A = {\log_2T\log_3T\over2V}$.

\medskip
{\bf Case 3}. When $\hf\t \log_2T\log_3T \le V \le (3\log 2T)/(8\log_22T)$
we take $A = 2/\t$.

\smallskip
Note that the last bound for $V$ comes from the bound (1.3) with $C=3/8$
(under the RH). Suppose that $\log |\zt| \ge V \ge 10\sqrt{\log_2T}$ holds.
Then Lemma 1 and Lemma 2 yield
$$
V \le S_1(t) + S_2(t) + {1+\lambda_0\over2A\t}V + O(\log_3T),\leqno(4.2)
$$
where we set
$$
S_1(t) := \Bigl|\sum_{p\le z}{\log(x/p)\over\log x}p^{-{1\over2}-
{\lambda_0\over\log x}-it} \Bigr|,\;
S_2(t) := \Bigl|\sum_{z<p\le x}{\log(x/p)\over\log x}p^{-{1\over2}-
{\lambda_0\over\log x}-it} \Bigr|.\leqno(4.3)
$$
This means that either
$$
S_1(t) \;\ge\; V_1 = V\left(1 - {7\over8A\t}\right) \leqno(4.4)
$$
or
$$
S_2(t) \;\ge\; {V\over8A\t}.\leqno(4.5)
$$
Namely, if neither (4.4) nor (4.5) is true, then (4.2) implies
that for some constant $C>0$
$$
V \le V\left(1 - {7\over8A\t}\right) + {V\over8A\t} + C\log_3T.
$$
Therefore we should have
$$
{V\over A} \;\ll\; \log_3T,\leqno(4.6)
$$
but (4.6) obviously cannot hold in view of the ranges of $V$ and the choice
of $A$ in Cases 1-3. Let now $\mu_i(T,H,V)\;(i=1,2)$ denote the measure of the
set of points  $t\in [T,\,T+H]$ for which  (4.4) and (4.5) hold,
respectively. Supposing that (4.4) holds then, by using Lemma 3
with $a(p) = {\log(x/p)\over\log x}p^{-\lambda_0/\log x}$, we obtain
$$
\mu_1(T,H,V)V_1^{2r} \le \int_T^{T+H}|S_1(t)|^{2r}\d t \ll Hr!
{\Bigl(\,\sum_{p\le z}{1\over p}\,\Bigr)}^r. \leqno(4.7)
$$
The condition in Lemma 3 ($x^r \le H$ with $x = z$) is equivalent to
$$
{Ar\over V\log_2T} \;\le\;1. \leqno(4.8)
$$
Recalling that
$$
\sum_{p\le X}{1\over p} = \log_2 X + O(1),
$$
it follows that
$$
\log z = {\log x\over \log_2T} = {A\t\over V\log_2T}\log T \le {\log T\over \log_2T}.
$$
since $A \le V$ in all cases. Therefore we have
$$
\sum_{p\le z}{1\over p} \;\le\; \log_2T\qquad(T\ge T_0).\leqno(4.9)
$$
Noting that Stirling's formula yields $r! \ll r^r\sqrt{r}{\roman e}^{-r}$, we
infer from (4.7) and (4.9) that
$$
\mu_1(T,H,V) \;\ll\; H\sqrt{r}
{\Bigl({r\log_2T\over{\roman e}V_1^2}\Bigr)}^r.
\leqno(4.10)
$$
In the Cases 1. and 2. and also in the Case 3. when
$V \le {2\over\t}\log_2^2T$, one chooses
$$
r \= \left[{V_1^2\over \log_2T}\right] \;\Bigl(\ge1\Bigr).
$$
With this choice of $r$ it is readily seen that (4.8) is satisfied, and (4.10) gives
$$
\mu_1(T,H,V) \;\ll\; H{\sqrt{V}\over\log_2T}
\exp\Bigl(-{V_1^2\over\log_2T}\Bigr).\leqno(4.11)
$$
Finally in the Case 3. when ${2\over\t}\log_2^2T \le V \le
(3\log 2T)/(8\log_2 2T)$ and $A = 2/\t$, we have
$$
V_1 = V\left(1 - {7\over8A\t}\right) = V\left(1 - {7\over16}\right)
> {V\over2},
$$
so that with the choice $r = [V/2]$  we see that (4.8) is again satisfied and
$$
\sqrt{r}{\Bigl({r\log_2T\over{\roman e}V_1^2}\Bigr)}^r
\le \sqrt{V}\,\left({2\log_2T\over {\roman e}V}\right)^r \le
V^{{1\over2}-{r\over4}} \ll \exp(-{\txt{1\over10}}V\log V),
$$
giving in this case
$$
\mu_1(T,H,V) \;\ll\; H\exp(-{\txt{1\over10}}V\log V).\leqno(4.12)
$$
\medskip
We bound $\mu_2(T,H,V) $ in a similar way by using (4.5). It follows,
again by Lemma 3, that
$$\eqalign{
\Bigl({V\over8A\t}\Bigr)^{2r}\mu_2(T,H,V) &\le
\int_T^{T+H}|S_2(t)|^{2r}\d t\cr&
\ll Hr!{\Bigl(\,\sum_{z<p\le x}{1\over p}\,\Bigr)}^r
= Hr!{\bigl(\log_2x - \log_2z + O(1)\bigr)}^r\cr&
\ll H\Bigl(r(\log_3T + O(1))\Bigr)^r.\cr}
$$
We obtain
$$
\mu_2(T,H,V) \ll H\Bigl({8A\over V}\Bigr)^{2r}\bigl(2r\log_3T\bigr)^r
\ll H\exp\left(-{V\over 2A}\,\log V\right).
\leqno(4.13)
$$
Namely the second inequality in (4.13) is equivalent to
$$
\left({A\over V}\right)^2\,r\log_3T \;\ll\; \exp\left(-{V\over2rA}\log V\right).
\leqno(4.14)
$$
In all Cases 1.-3. we take
$$
r = \left[{V\over A} - 1\right] \;\Bigl(\ge1\Bigr).
$$
The condition $x^r \le H$ in Lemma 3 is equivalent to $rA \le V$, which is
trivial with the above choice of $r$. To establish (4.14) note first that
$$
\left({A\over V}\right)^2\,r\log_3T \le {A\over V}\log_3T.\leqno(4.15)
$$
In the Case 1. the second expression in (4.15) equals $\log_3^2T/(2V)$, while
$$
\exp\left(-{V\over2rA}\log V\right) = \exp\left(-(\hf + o(1))\log V\right)
= V^{-1/2+o(1)}.
$$
Therefore it suffices to have
$$
{\log_3^2T\over V} \ll V^{-1/2+o(1)},
$$
which is true since $10\sqrt{\log_2T} \le V$. In the Case 2. the analysis is similar.
In the Case 3. we have $A = 2/\t$, hence $(A\log_3T)/V \ll (\log_3T)/V$ and
$$
\exp\left(-{V\over2rA}\log V\right) = \exp\Bigl(-\bigl({1\over2}\t+o(1)\bigr)\log V\Bigr)
= V^{-(\hf\t+o(1))},
$$
so that (4.14) follows again. Thus we have shown that in all cases
$$
\mu_2(T,H,V) \ll H\exp\left(-{V\over 2A}\,\log V\right).
\leqno(4.16)
$$
Theorem 2 follows now from (4.10), (4.11) and (4.16). Namely in the Case 1. we
have
$$
{V_1^2\over\log_2T} = {V^2\Bigl(1 - {7\over4\t\log_3T}\Bigr)^2
\over\log_2T} \le V\,{\log V\over\log_3T} = {V\log V\over2A},
$$
which gives (1.7). If the Case 2. holds we have again
$$
{V_1^2\over\log_2T} = {V^2\Bigl(1 - {7V\over4\t\log_2T\log_3T}\Bigr)^2\over\log_2T}
\le {V^2\log V\over \log_2T\log_3T} = {V\log V\over 2A},
$$
and (1.8) follows. In the Case 3. when
$\hf\t\log_2T\log_3T \le V \le {2\over\t}\log_2^2T$ we have
$$\eqalign{
\mu(T,H,V) &\ll H\exp\Bigl(-{V_1^2\over\log_2T}\Bigr) + H\exp(-\t V\log V)\cr&
\ll H\exp(-{\t\over20}V\log V),\cr}
$$
since
$$
{V_1^2\over\log_2T} \ge {V^2\over4\log_2T} \ge {\t V\log_2T\log_3T\over8\log_2T}
\ge {\t\over20}\,V\log V.
$$
In the remaining range of Case 3. we have
$$
\eqalign{
\mu(T,H,V) &\ll H\exp(-\txt{1\over10}V\log V) +
H\exp\bigl(-{V\over 2A}\,\log V\bigr)\cr&
\ll H\exp(-\txt{1\over10}V\log V) + \exp(-\txt{\t\over4}V\log V),\cr}
$$
and (1.9) follows. The proof of Theorem 2 is complete.


\vfill
\eject
\topskip1cm

\Refs

\item{[1]} F.V. Atkinson, The mean value of the Riemann zeta-function,
Acta Math. {\bf81}(1949), 353-376.

\item{[2]} J.B. Conrey, D.W. Farmer, J.P. Keating, M.O. Rubinstein
and N.C. Snaith, Integral moments of $L$-functions, Proc. Lond. Math.
Soc., III. Ser. {\bf91}(2005), 33-104.

\item{[3]} A. Ivi\'c, The Riemann zeta-function, John Wiley \&
Sons, New York 1985 (2nd edition. Dover, Mineola, New York, 2003).

\item{[4]} A. Ivi\'c,  The mean values of the Riemann zeta-function,
LNs {\bf 82}, Tata Inst. of Fundamental Research, Bombay 1991 (distr. by
Springer Verlag, Berlin etc.).

\item{[5]} M.B. Milinovich, Upper bounds for the moments of $\z'(\rho)$,
to appear, see {\tt arXiv:0806.0786}.

\item {[6]}  Y. Motohashi,  Spectral theory of the Riemann
zeta-function,  Cambridge University Press, Cambridge, 1997.

\item{[7]} K. Ramachandra, On the mean-value and omega-theorems
for the Riemann zeta-function, Tata Inst. Fundamental Research,
distr. by Springer Verlag), Bombay, 1995.

\item{[8]} K. Soundararajan, Moments of the Riemann zeta-function,
Annals of Math. {\bf170}(2009),981-993.

\item{[9]} E.C. Titchmarsh, The theory of the Riemann
zeta-function (2nd edition),  University Press, Oxford, 1986.

\vskip2cm
\endRefs

\enddocument

\bye